\documentclass[10pt,reqno,a4paper]{amsart}
\usepackage{amsfonts,amsmath,times}
\usepackage{tikz}

\allowdisplaybreaks

\newtheorem{theorem}{Theorem}

\DeclareMathOperator{\Betarv}{Beta}

\title{Analysis of some exactly solvable diminishing urn models}

\author[H.-K.~Hwang]{Hsien-Kuei Hwang}
\address{Hsien-Kuei Hwang\\
Institute of Statistical Science\\
Academia Sinica\\
Taipei, 115\\
Taiwan}
\email{hkhwang@stat.sinica.edu.tw}

\author[M.~Kuba]{Markus Kuba}
\address{Markus Kuba\\
Institut f{\"u}r Diskrete Mathematik und Geometrie\\
Technische Universit\"at Wien\\
Wiedner Hauptstr. 8-10/104\\
1040 Wien, Austria} %
\email{markus.kuba@tuwien.ac.at}

\author[A.~Panholzer]{Alois Panholzer}
\address{Alois Panholzer\\
Institut f{\"u}r Diskrete Mathematik und Geometrie\\
Technische Universit\"at Wien\\
Wiedner Hauptstr. 8-10/104\\
1040 Wien, Austria}
\email{Alois.Panholzer@tuwien.ac.at}

\thanks{This work was partially supported by the Austrian Science
Foundation FWF, grant S9608-N13.}

\keywords{Urn models, Generating functions, Limiting distribution}%
\subjclass[2000]{05A15,60F05,05C05} %

\begin{document}

\begin{abstract}
We study several exactly solvable P{\'{o}}lya-Eggenberger urn models
with a \emph{diminishing} character, namely, balls of a specified
color, say $x$ are completely drawn after a finite number of draws.
The main quantity of interest here is the number of balls left when
balls of color $x$ are completely removed. We consider several
diminishing urns studied previously in the literature such as the
pills problem, the cannibal urns and the OK Corral problem, and
derive exact and limiting distributions. Our approach is based on
solving recurrences via generating functions and partial
differential equations.
\end{abstract}

\maketitle

\vspace*{-3ex}

\begin{quote}
\renewcommand{\baselinestretch}{0.8}\small\normalsize
\begin{footnotesize}
On se propose d'\'etudier plusieurs mod{\`{e}}les d'urnes de
P{\'{o}}lya-Eggenberger de nature ``diminuante'' ayant des solutions
exactes, c'est-\`a-dire, les boules de couleur, disons $x$, sont
toutes prises apr\`es un nombre fini de tir\'ees. La quantit\'e\
principale qui nous interesse est le nombre de boules qui restent
dans l'urne au moment o\`u\ il n'y a plus de boules de couleur $x$.
Nous traitons, en particulier, plusieurs mod\`eles d'urnes
diminuantes propos\'es dans la litterature, comme le probl\`eme de
pillules, le mod\`ele d'urnes dit ``cannibaliste" et le probl\`eme
d'OK Corral, et obtenons des r\'esultats exactes et asymptotiques.
L'approche que nous utilisons est fond\'ee sur le traitement de
r\'ecurrences par voie de fonctions g\'en\'eratrices et \'equations
aux d\'erivativ\'es partielles.

\renewcommand{\baselinestretch}{1.0}\small\normalsize
\end{footnotesize}
\end{quote}

\section{Introduction}

\subsection{Diminishing urn models}

We are concerned here with the so-called P{\'{o}}lya-Eggenberger urn
models, which in the simplest case of two types of colors for the
balls  can be described as follows. At the beginning, the urn
contains $m$ black and $n$ white balls. At every step, we choose a
ball at random from the urn, examine its color and put it back into
the urn and then add/remove balls according to its color by the
following rules. If the ball is white, then we put $a$ white and $b$
black balls into the urn, while if the ball is black, then $c$ white
balls and $d$ black balls are put into the urn. The values $a, b, c,
d \in \mathbb{Z}$ are fixed integer values and the urn model is
specified by the transition matrix $M = \bigl(\begin{smallmatrix} a & b \\
c & d\end{smallmatrix}\bigr)$. Urn models with $r$ ($\ge2$) types of
colors can be described in an analogous way and are specified by an
$r \times r$ transition matrix.

Urn models are simple, useful mathematical tools for describing many
evolutionary processes in diverse fields of application such as
analysis of algorithms and data structures, statistics and genetics.
Due to their importance in applications, there is a huge literature
on the stochastic behavior of urn models; see for example
\cite{JohKot1977,KotBal1997}. Recently, a few different approaches
have been proposed, which yield deep and far-reaching results for
very general urn models; see
\cite{FlaDumPuy2006,FlaGabPek2005,Jan2004}.

Most papers in the literature impose the so-called \emph{tenability}
condition on the transition matrix, so that the process can be
continued \emph{ad infinitum} (or no balls of a given color being
completely removed). However, in some applications (examples given
below), there are urn models with a very different nature, which we
will refer to as ``diminishing urn models.'' For simplicity of
presentation, we describe them in the case of balls with two types
of colors, black and white. We consider P{\'{o}}lya-Eggenberger urn
models specified by a transition matrix $M =
\bigl(\begin{smallmatrix} a & b \\ c & d\end{smallmatrix}\bigr)$,
and in addition there is a set of absorbing states $\mathcal{S}
\subseteq \mathbb{N} \times \mathbb{N}$. The urn contains $m$ black
balls and $n$ white balls at the beginning and evolves by successive
draws at discrete instance according to the transition matrix until
an absorbing state $s = (j,k) \in \mathcal{S}$ is reached, namely,
the urn contains exactly $j$ black balls and $k$ white balls. Then
the urn process stops. We only call an urn model ``diminishing urn
model'' if it is guaranteed that from any initial state $(m,n) \in
\mathbb{N} \times \mathbb{N}$ (starting with $m$ black balls and $n$
white balls) we will reach an absorbing state $s \in \mathcal{S}$
after a finite number of draws.

Diminishing urn models with more than two type of balls can be
considered similarly; an example will be given below. For
diminishing urns, the main questions are (i) starting at state
$(m,n)$, what is the probability of reaching the absorbing state
$(j,k) \in \mathcal{S}$?, and (ii) what is then the number of balls
left?.

Motivated by concrete applications, we distinguish the following two
types of urns.
\begin{description}

\item[\textbf{Type A}] The entries of $M$ satisfy $a, b \le 0$,
$(a,b)\neq(0,0)$, $d<0$ and $c>0$, and the set of absorbing states
$\mathcal{S}$ consists of the vertical axis $m=0$ (or a vertical
wall $0 \le m \le M$, $M \ge 0$).

\item[\textbf{Type B}] The entries of $M$ satisfy $a, b, c, d
\le 0$, $(a,b) \neq (0,0)$ and $(c,d)\neq (0,0)$ and the set of
absorbing states $\mathcal{S}$ consists of the vertical axis $m = 0$
(or a vertical wall $0 \le m \le M$, $M \ge 0$) and the horizontal
axis $n=0$ (or a horizontal wall $0 \le n \le N$, $N \ge 0$).

\end{description}
Note that the conditions on Type A urn models are in general not
sufficient to guarantee that an absorbing state will be reached
(if $b<-1$ then the urn process could reach states with $n<0$), but
this is the case for all models we consider here.

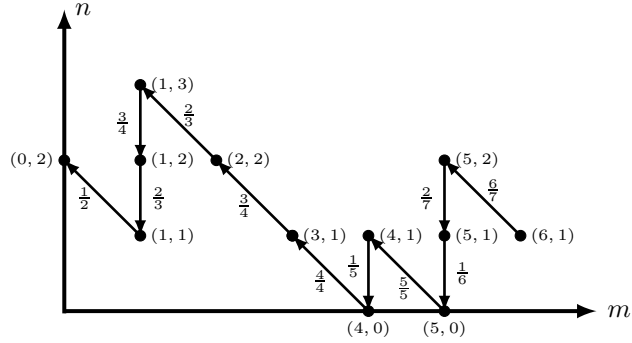
\begin{figure}
\centering
\begin{tikzpicture}
\draw[line width=1.5pt,-latex]
    (0,0) -- (7,0) node[right] {$m$};%
\draw[line width=1.5pt,-latex]
    (0,-0.03) -- (0,4) node[pos=1,right] {$n$};%
\draw[line width=1pt,-latex] (6,1) -- (5,2)
    node[right,pos=0] {\tiny $(6,1)$}%
    node[right,pos=.6] {\tiny $\frac67$};%
\draw[line width=1pt,-latex] (5,2) -- (5,1)
    node[right,pos=0] {\tiny $(5,2)$}%
    node[left,pos=.5] {\tiny $\frac27$};%
\draw[line width=1pt,-latex] (5,1) -- (5,0)
    node[right,pos=0] {\tiny $(5,1)$}%
    node[right,pos=.5] {\tiny $\frac16$};%
\draw[line width=1pt,-latex] (5,0) -- (4,1)
    node[below,pos=0] {\tiny $(5,0)$}%
    node[left,pos=.3] {\tiny $\frac55$};%
\draw[line width=1pt,-latex] (4,1) -- (4,0)
    node[right,pos=0] {\tiny $(4,1)$};%
\draw (3.8,0.6) node {\tiny $\frac15$};%
\draw[line width=1pt,-latex] (4,0) -- (3,1)
    node[below,pos=0] {\tiny $(4,0)$}%
    node[left,pos=.4] {\tiny $\frac44$};%
\draw[line width=1pt,-latex] (3,1) -- (2,2)
    node[right,pos=0] {\tiny $(3,1)$}%
    node[left,pos=.4] {\tiny $\frac34$};%
\draw[line width=1pt,-latex] (2,2) -- (1,3)
    node[right,pos=0] {\tiny $(2,2)$}%
    node[right,pos=.6] {\tiny $\frac23$};%
\draw[line width=1pt,-latex] (1,3) -- (1,2)
    node[right,pos=0] {\tiny $(1,3)$}%
    node[left,pos=.5] {\tiny $\frac34$};%
\draw[line width=1pt,-latex] (1,2) -- (1,1)
    node[right,pos=0] {\tiny $(1,2)$}%
    node[right,pos=.5] {\tiny $\frac23$};%
\draw[line width=1pt,-latex] (1,1) -- (0,2)
    node[right,pos=0] {\tiny $(1,1)$}%
    node[left,pos=.5] {\tiny $\frac12$}%
    node[left,pos=1] {\tiny $(0,2)$};%
\draw[fill] (6,1) circle (2pt);%
\draw[fill] (5,2) circle (2pt);%
\draw[fill] (5,1) circle (2pt);%
\draw[fill] (5,0) circle (2pt);%
\draw[fill] (4,1) circle (2pt);%
\draw[fill] (3,1) circle (2pt);%
\draw[fill] (4,0) circle (2pt);%
\draw[fill] (3,1) circle (2pt);%
\draw[fill] (2,2) circle (2pt);%
\draw[fill] (1,3) circle (2pt);%
\draw[fill] (1,2) circle (2pt);%
\draw[fill] (1,1) circle (2pt);%
\draw[fill] (0,2) circle (2pt);%
\end{tikzpicture}
\caption{An example of a weighted path from $(6,1)$ to the absorbing
state $(0,2)$ for the so called pills problem with transition matrix
$M=[-1, 0; 1, -1]$ and the vertical absorbing axis $\mathcal{S} =
\{(0,n): n \ge 0\}$. The illustrated path has weight $\frac{6}{7}
\frac{2}{7} \frac{1}{6} \frac{5}{5} \frac{1}{5} \frac{4}{4}
\frac{3}{4} \frac{2}{4} \frac{3}{4} \frac{2}{3} \frac{1}{2} =
\frac{3}{3920}$.\label{fig1}}
\end{figure}

It is helpful to describe the evolution of the urn model by weighted
lattice paths, which is described in the case of urns with two types
of balls. If the urn contains $m$ black balls and $n$ white balls
and we select a white ball (with probability $\frac{n}{m+n}$), then
this corresponds to a step from $(m,n)$ to $(m+a,n+b)$, to which the
weight $\frac{n}{m+n}$ is associated; and if we select a black ball
(with probability $\frac{m}{m+n}$), this corresponds to a step from
$(m,n)$ to $(m+c,n+d)$ (with weight $\frac{m}{m+n}$). The weight of
a path after $t$ successive draws consists of the product of the
weight of every step. By this correspondence, the probability of
starting at $(m,n)$ and ending at $(j,k)$ is equal to the sum of the
weights of all possible paths starting at state $(m,n)$ and ending
at the absorbing state $(j,k) \in \mathcal{S}$ (which did not reach
any absorbing state before). Unfortunately, the expressions so
obtained for the probability are, although exact, less useful for
large $m$ or $n$. An example for the weighted path corresponding to
the evolution of a diminishing urn is given in Figure~\ref{fig1}.

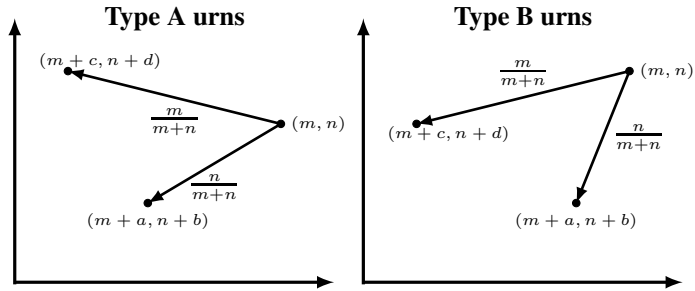
\begin{figure}
\centering
\begin{tikzpicture}[scale=0.7]
\draw[line width=1.2pt,-latex] (0,0) -- (6,0);%
\draw[line width=1.2pt,-latex] (0,-0.04) -- (0,5);%
\draw[line width=1pt,-latex] (5,3) -- (1,4)
node[above,pos=.85] {\tiny $(m+c,n+d)$}%
node[right,pos=0] {\tiny $(m,n)$}%
node[below,pos=0.5] {\small $\frac{m}{m+n}$};%
\draw[line width=1pt,-latex] (5,3) -- (2.5,1.5)
node[below,pos=1] {\tiny $(m+a,n+b)$}%
node[below,pos=0.5] {\small $\frac{n}{m+n}$};%
\draw[fill] (5,3) circle (2pt);%
\draw[fill] (1,4) circle (2pt);%
\draw[fill] (2.5,1.5) circle (2pt);%
\draw (3,5) node {\textbf{Type A urns}};
\end{tikzpicture}
\begin{tikzpicture}[scale=0.7]
\draw[line width=1.2pt,-latex] (0,0) -- (6,0);%
\draw[line width=1.2pt,-latex] (0,-0.04) -- (0,5);%
\draw[line width=1pt,-latex] (5,4) -- (1,3)
node[below,pos=.85] {\tiny $(m+c,n+d)$}%
node[right,pos=0] {\tiny $(m,n)$}%
node[above,pos=0.5] {\small $\frac{m}{m+n}$};%
\draw[line width=1pt,-latex] (5,4) -- (4,1.5)
node[below,pos=1] {\tiny $(m+a,n+b)$}%
node[right,pos=0.5] {\small $\frac{n}{m+n}$};%
\draw[fill] (5,4) circle (2pt);%
\draw[fill] (1,3) circle (2pt);%
\draw[fill] (4,1.5) circle (2pt);%
\draw (3,5) node {\textbf{Type B urns}};
\end{tikzpicture}
\caption{Type A and Type B urn models.\label{fig2}}
\end{figure}

The description of the urn model via weighted lattice paths then
gives the following interpretation of the Type A and Type B urn
models: in Type A we have a step lying in the lower-left quadrant
and one step lying in the upper-left quadrant, whereas in Type B
both steps are lying in the lower-left quadrant; see
Figure~\ref{fig2}.

\subsection{Examples\label{ssec12}}

We first describe a few motivating examples of diminishing urn
models.

\paragraph{\emph{\textbf{The pills problem.}}} The transition matrix
is given by $M = \bigl(\begin{smallmatrix} -1 & 0 \\ 1 &
-1\end{smallmatrix}\bigr)$ and the absorbing axis is $\mathcal{S} =
\{(0,n): n \ge 0\}$. An interpretation is as follows. An urn has two
types of pills in it, which are single-unit and double-unit pills,
respectively. At every step, we pick a pill uniformly at random. If
a single-unit pill is chosen, then we eat it up, and if the pill is
of double unit, we break it into two halves---one half is eaten up
and the other half is now considered of single unit and thrown back
into the urn. The question is then, when starting with $n$
single-unit pills and $m$ double-unit pills, what is the probability
that $k$ single-unit pills remain in the urn when all double-unit
pills are drawn?

This problem has been stated in \cite{KnuMcC1991}, where the authors
asked for a formula for the expected number of remaining single-unit
pills, when there are no double-unit pills in the urn. The solution
appeared in \cite{Hes1992}. A more refined study is given recently
in \cite{BrePro2003}, where they derive exact formul{\ae} for the
variance and the third moment of the number of remaining single-unit
pills; furthermore, a few generalizations are proposed.

A natural generalization is to consider $r$ types of pills, which
are of $i$ units, $i=1,\dots, r$, respectively. At every time step,
a pill is chosen uniformly at random; if the pill is of single unit,
it is eaten up, and if the pill is of $i$ units, $i \ge 2$, it is
broken into two parts, one of single unit and the other of $(i-1)$
units. The piece of single unit is eaten up and the remaining piece
is thrown back into the urn. We stop if there are no more pills of
the largest units ($r$).

This problem corresponds to the diminishing urn model with the $r
\times r$-transition matrix
\begin{small}
\begin{equation*}
    M = \left(
    \begin{smallmatrix}
        -1 & 0 & 0 & \cdots & 0 & 0 & 0 \\[-1ex]
        1 & -1 & 0 & \ddots & \ddots & \ddots & 0 \\[-1ex]
        0 & 1 & -1 & \ddots & \ddots & \ddots & 0 \\[-1ex]
        \vdots & \ddots & \ddots & \ddots & \ddots &
        \ddots & \vdots \\[-1ex]
        0 & \ddots & \ddots & \ddots & -1 & 0 & 0 \\[-1ex]
        0 & \ddots & \ddots & \ddots & 1 & -1 & 0 \\
        0 & 0 & 0 & \cdots & 0 & 1 & -1
    \end{smallmatrix}
    \right)
\end{equation*}
\end{small}
and the absorbing hyperplane $\mathcal{S} = \{(n_{1}, \dots,
n_{r-1}, 0): n_{1}, \dots, n_{r-1} \ge 0\}$. We will be interested
in finding the probability that $k$ pills of single unit remain in
the urn when there are no more pills of $r$ units, the starting
configuration being $n_i$ pills of $i$ units.

\paragraph{\emph{\textbf{A variant of the pills problem}}.}
To illustrate how a minor change in the entries of the transition
matrix leads to very different behavior, we will also consider the
transition matrix $M =\bigl(\begin{smallmatrix} -1 & 0 \\ 1 &
-2\end{smallmatrix}\bigr)$ and the absorbing wall $\mathcal{S} =
\{(0,n): n \ge 0\} \cup \{(1,n): n \ge 0\}$.

\paragraph{\emph{\textbf{The cannibal urn.}}}
Introduced by R.~F.~Greene (unpublished) and analyzed in details by
Pittel in \cite{Pit1987}, this urn model is a slight modification of
the diminishing urn with $M=\bigl(\begin{smallmatrix} 0 & -1
\\ 1 & -2\end{smallmatrix}\bigr)$ and the vertical wall of
absorbing states $\mathcal{S} = \{(0,n): n \ge 0\} \cup \{(1,n): n
\ge 0\}$. In terms of weighted lattice paths, one starts at position
$(m,n)$, the weight (and thus the probability) of a step to
$(m-1,n)$ is $\frac{n}{m-1+n}$ (not $\frac{n}{m+n}$), and the weight
to $(m-2,n+1)$ is $\frac{m-1}{m-1+n}$. The approach we use is also
applicable to this modified urn model.

Such an urn was introduced to model the behavior of cannibals in
biological population. It can be described as follows. A population
consists of cannibals and non-cannibals. At every time step, a
non-cannibal is selected as victim and removed; after that a member
in the remaining population (cannibals and non-cannibals) is
selected uniformly at random. If the selected individual is a
cannibal it remains as a cannibal, but if the selected individual is
a non-cannibal, it becomes then a cannibal. The question is, when
starting with $n$ cannibals and $m$ non-cannibals, what is the
number of resulting cannibals in the population at the moment when
all non-cannibals are removed?

\paragraph{\emph{\textbf{The OK Corral problem.}}} This corresponds
to the urn $M =\bigl(\begin{smallmatrix} 0 & -1 \\ -1 &
0\end{smallmatrix}\bigr)$ with two absorbing axes: $\mathcal{S} =
\{(0,n): n \ge 0\} \cup \{(m,0): m \ge 0\}$. An interpretation is as
follows. Two groups of gunmen, group A and group B (with $n$ and $m$
gunmen, respectively), face each other. At every discrete time step,
one gunman is chosen uniformly at random who then shoots and kills
exactly one gunman of the other group. The bloody gunfight ends when
one group gets completely ``eliminated". Two questions are of
interest: (i) what is the probability that group A (group B)
survives? and (ii) what is the probability that the gunfight ends
with $k$ survivors of group A (group B)?

This problem was introduced by Williams and McIlroy in
\cite{WilMcI1998} and studied recently by several authors using
different approaches, leading to very interesting results; see
\cite{FlaDumPuy2006,Kin1999,KinVol2003}. Also the urn corresponding
to the OK corral problem can be viewed as a basic model in the
mathematical theory of warfare and conflicts; see \cite{KinVol2003}.

\paragraph{\emph{\textbf{Sampling without replacement.}}} This
is a toy example and corresponds to the urn $M =
\bigl(\begin{smallmatrix} -1 & 0 \\ 0 & -1\end{smallmatrix}\bigr)$
with two absorbing axes: $\mathcal{S} = \{(0,n): n \ge 0\} \cup
\{(m,0): m \ge 0\}$. In this classical model, balls are drawn one
after another from an urn containing balls of two different colors
and not replaced. What is the probability that $k$ balls of one
color remain when balls of the other color are all removed?

\subsection{Recurrence\label{ssec13}}

For diminishing urns, we study the position of the absorbing state.
Probabilistically, we consider the pair of random variables
$(X_{n,m}^{(1)},X_{n,m}^{(2)})$, such that
$\mathbb{P}\{(X_{n,m}^{(1)},X_{n,m}^{(2)}) = (j,k)\}$ gives the
probability that when starting at state $(m,n)$ (with $m$ black
balls and $n$ white balls), the urn process reaches the absorbing
state $(j,k)$, namely, the process terminates with $j$ black balls
and $k$ white balls. For diminishing urns with a single vertical
absorbing axis (or wall), we are only interested in the vertical
position of the absorbing state; so we define $X_{n,m} :=
X_{n,m}^{(2)}$ and $\mathbb{P}\{X_{n,m} = k\}$ is then the
probability that when starting with $m$ black balls and $n$ white
balls, the urn process stops with $k$ white balls remaining in the
urn. We consider the probability generating function $h_{n,m}(v_{1},
v_{2})$ or $h_{n,m}(v)$, respectively, defined by
\begin{equation}
\label{eqna2}
    h_{n,m}(v_{1},v_{2}) := \sum_{j \ge 0} \sum_{k \ge 0}
    \mathbb{P}\big\{(X_{n,m}^{(1)},X_{n,m}^{(2)}) = (j,k)\big\}
    v_{1}^{j} v_{2}^{k}, \quad h_{n,m}(v) := \sum_{k \ge 0}
    \mathbb{P}\{X_{n,m}=k\} v^{k}.
\end{equation}
According to the outcome of the first draw of the urn process, we
obtain the following recurrences for the probability generating
functions
\begin{subequations}
\label{eqna1}
\begin{align}
    h_{n,m}(v_{1},v_{2}) & =
    \frac{n}{m+n} h_{n+a,m+b}(v_{1},v_{2})
    + \frac{m}{m+n} h_{n+c,m+d}(v_{1},v_{2}), \label{eqna1a} \\
    h_{n,m}(v) & = \frac{n}{m+n} h_{n+a,m+b}(v)
    + \frac{m}{m+n} h_{n+c,m+d}(v), \label{eqna1b}
\end{align}
\end{subequations}
for $(m,n) \not\in S$. The boundary values at the absorbing states
$(m,n) \in S$ are given by $h_{n,m}(v_{1},v_{2}) = v_{1}^{m}
v_{2}^{n}$ and $h_{n,m}(v) = v^{n}$, respectively.

We solve such recurrences via generating functions for a few special
cases below. For urn models of Type B we can always introduce
generating functions\footnote{The generating function also depends
on $v_{1}$ and $v_{2}$, but we avoid the heavier notation
$H(z,w;v_{1},v_{2}$).}
\begin{equation*}
    H(z,w) := \sum_{(m,n) \not\in S}
    h_{n,m}(v_{1},v_{2}) z^{n} w^{m},
\end{equation*}
and the recurrence~\eqref{eqna1a} can be translated into the
following first order linear partial differential equation (PDE)
\begin{equation*}
    z (1-z^{-a}w^{-b}) H_{z}(z,w)
    + w(1-z^{-c}w^{-d})H_{w}(z,w)
    + (a z^{-a}w^{-b}+dz^{-c}w^{-d}) H(z,w) = F(z,w),
\end{equation*}
(see \cite{FlaGabPek2005}) where the inhomogeneous part $F(z,w)$ is
fully determined by the boundary values. Such PDEs can be treated
(at least in principle) by the method of characteristics, see, for
example, \cite{Tay1996}. For urn models of Type A, the situation
becomes more involved. The same approach may still apply but the
additional difficulty is the fact that the inhomogeneous part
$F(z,w)$ involves evaluations of the function $H(z,w)$ or their
partial derivatives at $z=0$. Fortunately, for the cases we consider
here, we can solve this problem by introducing an appropriate
normalizing factor; the resulting generating function satisfies then
a simpler PDE (with boundary values properly eliminated) that can be
explicitly solved.

Another general difficulty in solving the recurrences~\eqref{eqna1}
by solving the associated PDEs is how to adapt the general solution
to the boundary values. By the method of characteristics, we see
that the general solution is given by $H(z,w) = H^{[p]}(z,w) +
f(z,w) C(\xi(z,w))$ with an arbitrary continuous function $C(x)$.
Often it is not obvious how to find $C(x)$ such that $H(z,w)$
satisfies the boundary values. However, for the examples treated
here, we can always solve this problem by using the analyticity of
the function $H(z,w)$ in a neighborhood of $(z,w)=(0,0)$, by
choosing a suitable curve $z=q(w)$ and by considering the limit
$\lim_{w \to 0}H(q(w),w)$ (depending on $f(z,w)$ and $\xi(z,w)$).

As we show later, we obtain for all problems mentioned above
closed-form solutions for $H(z,w)$. From such exact forms, we can
easily derive the corresponding exact solutions for the underlying
probability. Also we can apply general analytic tools such as
singularity analysis and saddle-point method
(see~\cite{FlaSed2006+}) and obtain rather precise information on
the asymptotic growth of the underlying probabilities. However, for
problems in two or more variables as we are dealing with here, the
treatment is generally more involved than in the univariate case.

\section{The pills problem}

\subsection{The original problem}

We start by considering Type A diminishing urn model with the
transition matrix $M = \bigl(\begin{smallmatrix} -1 & 0 \\ 1 &
-1\end{smallmatrix}\bigr)$ and the vertical absorbing axis
$\mathcal{S} = \{(0,n): n \ge 0\}$.

The recurrence \eqref{eqna1b} for the probability generating
function $h_{n,m}(v)$ now becomes
\begin{equation}
\label{eqnd1}
    h_{n,m}(v) = \frac{n}{n+m} h_{n-1,m}(v) + \frac{m}{n+m}
    h_{n+1,m-1}(v),
\end{equation}
for $n \ge 0$ and $m \ge 1$, with the boundary values $h_{n,0}(v) =
v^{n}$.

Instead of considering the generating function $\tilde{H}(z,w) :=
\sum_{n \ge 0} \sum_{m \ge 1} h_{n,m}(v) z^{n} w^{m}$, which will
involve the unknown boundary values $h_{0,m}(v)$ (or
$\tilde{H}(0,w)$) in the resulting PDE, we introduce the modified
generating function
\begin{equation}
\label{eqnd2}
    H(z,w) := \sum_{n \ge 0} \sum_{m \ge 1}
    \binom{m+n}{m} h_{n,m}(v) z^{n} w^{m}.
\end{equation}
Then $H$ satisfies, by recurrence~\eqref{eqnd1}, the first-order
linear PDE
\begin{equation}
\label{eqnd3}
    (z-z^{2}-w) H_{z}(z,w) + w(1-z) H_{w}(z,w) - z H(z,w)
    = \frac{wv}{(1-vz)^{2}},
\end{equation}
with the initial condition $H(z,0)=0$. We see that the unknown
boundary values $h_{0,m}(v)$ nicely disappear.

To solve equation \eqref{eqnd3}, we apply the method of
characteristics. Thus we first consider the corresponding reduced
PDE
\begin{equation}
\label{eqnd31}
    (z-z^{2}-w) H_{z}(z,w) + w(1-z) H_{w}(z,w) = 0,
\end{equation}
and find the first integrals for the system of ordinary differential
equations (the so-called system of characteristic differential
equations)
\begin{equation}
\label{eqnd4}
    \dot{z} = z - z^{2} - w, \quad \dot{w} = w (1-z).
\end{equation}
We regard here $z$ and $w$ as dependent variables of $t$, namely,
$z=z(t)$, $w=w(t)$ and $\dot{z} := \frac{d z(t)}{d t}$, etc. By
reducing \eqref{eqnd4} to a differential equation (DE) of Bernoulli
type, we obtain the following first integral of \eqref{eqnd4}
\begin{equation*}
    \xi(z,w) := \frac{w e^{z/w}}{1-z-w} = \text{const.}
\end{equation*}
Thus the general solution of the reduced PDE~\eqref{eqnd31} is as
follows.
\begin{equation*}
    H^{[r]}(z,w) = C\Big(\frac{w e^{z/w}}{1-z-w}\Big),
\end{equation*}
where $C(x)$ is an arbitrary continuous function

Now consider the inhomogeneous PDE
\begin{equation}
\label{eqnd5}
    (z-z^{2}-w) H_{z}(z,w) + w(1-z) H_{w}(z,w) - z H(z,w) = F(z,w).
\end{equation}
We use the following transformation from $(z,w)$-coordinates to
$(\eta, \xi)$-coordinates: $\xi = \frac{w e^{z/w}}{1-w-z}$ and $\eta
= \frac{w}{1-w-z}$, or equivalently $z = z(\eta,\xi) = \frac{\eta
\log(\xi/\eta)}{1+\eta+\eta\log(\xi/\eta)}$ and $w = w(\eta,\xi) =
\frac{\eta}{1+\eta+\eta\log(\xi/\eta)}$, which leads to the DE
\begin{equation}
\label{eqnd6}
    H_{\eta}(\eta, \xi) - \frac{\log(\xi/\eta)}
    {1+\eta+\eta \log(\xi/\eta)} H(\eta,\xi)
    = \frac{1}{\eta} F\big(z(\eta,\xi),w(\eta,\xi)\big).
\end{equation}
The general solution of the corresponding homogeneous DE
$H_{\eta}(\eta, \xi) - \log(\xi/\eta) H(\eta,\xi) /(1+\eta+\eta
\log(\xi/\eta)) = 0$ can be obtained easily and is given by
\begin{equation*}
    H^{[h]}(z,w) = \frac{1}{1-w-z} C\Big(\frac{w
    e^{\frac{z}{w}}}{1-w-z}\Big),
\end{equation*}
where we applied the inverse $(\eta,\xi)$-transform.

The inhomogeneous DE \eqref{eqnd6} can then be solved by using the
method of variation of parameters. We obtain for the inhomogeneous
part $F(z,w) =w v/(1-vz)^{2}$ the following particular solution
\begin{equation}
\label{eqnd7}
    H^{[p]}(z,w) = v w \int_{0}^{1}
    \frac{dq}{\big(1-z(1+(v-1)q)-w(1-q-(v-1)q \log q)\big)^{2}}.
\end{equation}
It turns out that the particular solution \eqref{eqnd7}, which is
analytic around $z=0$ and $w=0$, already satisfies the initial
condition, so \eqref{eqnd7} is the required solution of the problem,
$H(z,w) = H^{[p]}(z,w)$.

Extracting coefficients of $z^n$ and $w^m$ in \eqref{eqnd7} gives
for $n \ge 0$ and $m \ge 1$ the following explicit form
\begin{equation}
   \label{eqnd8}
   h_{n,m}(v) = \frac{1}{\binom{n+m}{n}} [z^{n} w^{m}] H(z,w) =
   m v \int_{0}^{1} (1+(v-1)q)^{n} (1-q-(v-1)q\log q)^{m-1} dq.
\end{equation}

From this expression, the expectation $\mathbb{E}(X_{n,m}) = h_{n,m}'(1)$ 
can be easily derived and is given by
\begin{equation}
   \label{eqnd32}
   \mathbb{E}(X_{n,m}) = \frac{n}{m+1} + H_{m};
\end{equation}
cf. \cite{BrePro2003,Hes1992}.

Higher moments can be obtained similarly by taking higher
derivatives from \eqref{eqnd8}, but the expressions soon become very
messy; see \cite{BrePro2003} for the second and the third moments.
Instead, we can apply \eqref{eqnd8} to derive the limiting
distribution of $X_{n,m}$, for all ranges of $n$ and $m$ satisfying
$\max(m,n) \to \infty$. The idea is roughly as follows. We first
compute asymptotic approximations to the $r$-th factorial moments
$\mathbb{E}(X_{n,m}^{\underline{r}}) := \mathbb{E}\big(X_{n,m}
(X_{n,m}-1) \cdots (X_{n,m}-r+1)\big)$ starting from the relation
$\mathbb{E}(X_{n,m}^{\underline{r}}) = h_{n,m}^{(r)}(1)$ and then 
by evaluating asymptotically the integrals
as derivatives of the Beta-function. The result is
\begin{align*}
    \mathbb{E}\big(X_{n,m}^{r}\big) \sim
    \mathbb{E}(X_{n,m}^{\underline{r}}) =
    \left\{\begin{array}{ll}  \displaystyle
    r!\Big(\frac{n}{m} + \log m\Big)^{r}
    \big(1+\mathcal{O}\big((\log m)^{-1}\big)\big),&
    \text{for} \; m \to \infty, \\[1.5ex]
    \displaystyle
    \frac{n^r}{\binom{m+r}{r}}
    \big(1+\mathcal{O}(n^{-1})\big), &
    \text{for} \; m \; \text{fixed} \;
    \text{and} \; n \to \infty.
    \end{array}\right.
\end{align*}
We then obtain the limiting distributions of $X_{n,m}$ after proper
normalization, justified by standard arguments (moment sequence
uniquely characterizes the distribution).

We collect our results for the pills problem in the following
theorem.
\begin{theorem}
Starting with $m$ double-unit pills and $n$ single-unit pills, the
probability generating function $h_{n,m}(v) := \sum_{k \ge 0}
\mathbb{P}\{X_{n,m}=k\} v^{k}$ of the number $X_{n,m}$ of the
remaining single-unit pills in the urn when all double-unit pills
are all taken is given by
\begin{equation*}
    h_{n,m}(v) = m v \int_{0}^{1} (1+(v-1)q)^{n}
    (1-q-(v-1)q\log q)^{m-1} dq.
\end{equation*}

If $m \to \infty$, then the random variable $X_{n,m}$ converges,
after suitable scaling, in distribution to an exponentially
distributed random variable $X$ with parameter $\lambda = 1$, namely
\begin{equation*}
    \frac{X_{n,m}}{\frac{n}{m}+\log m} \xrightarrow{(d)} X,
\end{equation*}
where $X$ has density $f(x) = e^{-x}$ for $x \ge 0$.

If $m$ is fixed and $n \to \infty$, then the random variable
$X_{n,m}$ converges, after suitable scaling, in distribution to a
Beta random variable $B_{m}$; in symbol
\begin{equation*}
    \frac{X_{n,m}}{n} \xrightarrow{(d)} B_{m} \stackrel{(d)}{=} \Betarv(1,m),
\end{equation*}
where $B_{m}$ has density $m (1-x)^{m-1}$, $0 \le x \le 1$.
\end{theorem}
Details of the proofs will be given in the full version of this
extended abstract.

\subsection{A generalization to $r$ pills}

We consider the random variable $X_{n_{1}, \dots, n_{r}}$, which
gives the number of single-unit pills when all pills of $r$ units
are all taken, starting with $n_i$ pills of $i$ units,
$i=1,\dots,r$. The probability generating function $h_{n_{1}, \dots,
n_{r}}(v) := \sum_{k \ge 0} \mathbb{P}\{X_{n_{1}, \dots, n_{r}} =
k\} v^{k}$ satisfies for $n_{1}, \dots, n_{r-1} \ge 0$, $n_{r} \ge
1$ the recurrence
\begin{equation}
\label{eqne1}
    h_{n_{1}, \dots, n_{r}}(v) =
    \frac{n_{1}}{n_1+\cdots+n_r}
    h_{n_{1}-1, n_{2}, \dots, n_{r}}(v)
    + \sum_{j=2}^{r} \frac{n_{j}}
    {n_1+\cdots+n_r} h_{n_{1}, \dots, n_{j-2}, n_{j-1}+1,
    n_{j}-1, n_{j+1}, \dots, n_{r}}(v),
\end{equation}
with the boundary value $h_{n_{1}, \dots, n_{r-1}, 0}(v) =
v^{n_{1}}$. Let
\begin{equation}
\label{eqne2}
    H(z_{1}, \dots, z_{r}) := \sum_{n_{1} \ge 0} \cdots
    \sum_{n_{r-1} \ge 0} \sum_{n_{r} \ge 1}
    \binom{n_{1}+\cdots+n_{r}}{n_{1}, \dots, n_{r}} h_{n_{1}, \dots,
    n_{r}}(v) z_{1}^{n_{1}} \cdots z_{r}^{n_{r}}.
\end{equation}
The recurrence~\eqref{eqne1} then translates into the first-order
linear PDE
\begin{align}
\label{eqne3}
    \sum_{j=1}^{r-1} & (z_{j}-z_{1} z_{j} - z_{j+1})
    H_{z_{j}}(z_{1}, \dots, z_{r}) + (z_{r}-z_{1} z_{r})
    H_{z_{r}}(z_{1}, \dots, z_{r}) - z_{1}
    H(z_{1}, \dots, z_{r}) \nonumber \\
    & = \frac{z_{r}}{(1-v z_{1}-z_{2}-\cdots-z_{r-1})^{2}},
\end{align}
for $r\ge3$, with the boundary condition $H(z_{1}, \dots, z_{r-1},
0) = 0$.

By the method of characteristics, we then consider the
characteristic system of DEs
\begin{equation}
\label{eqne4}
    \dot{z}_{1} = z_{1} - z_{1}^{2} - z_{2},
    \quad \dot{z}_{2} = z_{2} - z_{1} z_{2} - z_{3}, \quad
    \dots, \quad \dot{z}_{r-1} = z_{r-1} - z_{1} z_{r-1} - z_{r},
    \quad \dot{z}_{r} = z_{r} - z_{1} z_{r}.
\end{equation}
We can show that the $r-1$ functions $\xi_{1}(z_{1}, \dots, z_{r})$,
\dots, $\xi_{r-1}(z_{1}, \dots, z_{r})$ given below, where
$\xi_{1}$, \dots, $\xi_{r-2}$ are given implicitly as the solution
of a linear system of equations, give $r-1$ independent first
integrals of \eqref{eqne4}
\begin{small}
\begin{align*}
    \frac{z_{r-2}}{z_{r}} & =
    \frac{\big(\frac{z_{r-1}}{z_{r}}\big)^{2}}{2!} + \xi_{r-2}, \\
    \frac{z_{r-3}}{z_{r}} & =
    \frac{\big(\frac{z_{r-1}}{z_{r}}\big)^{3}}{3!} + \xi_{r-2}
    \frac{\big(\frac{z_{r-1}}{z_{r}}\big)}{1!} + \xi_{r-3}, \\
    \frac{z_{r-4}}{z_{r}} & =
    \frac{\big(\frac{z_{r-1}}{z_{r}}\big)^{4}}{4!} +
    \xi_{r-2} \frac{\big(\frac{z_{r-1}}{z_{r}}\big)^{2}}{2!}
    + \xi_{r-3} \frac{\big(\frac{z_{r-1}}{z_{r}}\big)}{1!}
    + \xi_{r-4}, \\
    \mbox{} \vdots \quad & = \quad\vdots \\
    \frac{z_{1}}{z_{r}} & =
    \frac{\big(\frac{z_{r-1}}{z_{r}}\big)^{r-1}}{(r-1)!} +
    \xi_{r-2} \frac{\big(\frac{z_{r-1}}{z_{r}}\big)^{r-3}}
    {(r-3)!} + \xi_{r-3}
    \frac{\big(\frac{z_{r-1}}{z_{r}}\big)^{r-4}}
    {(r-4)!}+ \cdots + \xi_{2}
    \frac{\big(\frac{z_{r-1}}{z_{r}}\big)}{1!} + \xi_{1}, \\
    \xi_{r-1} & = \frac{z_{r}}{1-z_{1}-\cdots-z_{r}}
    e^{\frac{z_{r-1}}{z_{r}}}.
\end{align*}
\end{small}
We can solve the PDE~\eqref{eqne3} by introducing $\eta :=
\frac{z_{r}}{1-z_{1}-\cdots-z_{r}}$ and $\xi_{1}$, \dots,
$\xi_{r-1}$ as above and applying a transform to the $(\eta,
\xi_{1}, \dots, \xi_{r-1})$-coordinates. We obtain then the
following explicit solution
\begin{multline}
   \label{eqna33}
    H(z_{1}, \dots, z_{r}) = \\
    z_{r} \int_{0}^{1}
    \frac{d q}{\Big(1-\sum_{j=1}^{r-1} (1-(-1)^{j-1} (1-v) q
    \frac{\log^{j-1} q}{(j-1)!}) z_{j}
    - (1-q-(-1)^{r-1}(1-v)q
    \frac{\log^{r-1}q}{(r-1)!})z_{r}\Big)^{2}}.
\end{multline}
Thus we obtain after extracting coefficients of \eqref{eqna33} an
explicit formula for $h_{n_{1}, \dots, n_{r}}(v)$, which is given by
the following theorem.
\begin{theorem}
Starting with $n_{1}$ pills of size $1$, \dots, $n_{r}$ pills of
size $r$, $r\ge3$, the probability generating function $h_{n_{1},
\dots, n_{r}}(v)$ of the number $X_{n_{1}, \dots, n_{r}}$ of pills
of single-unit pills remaining in the urn when all pills of $r$
units are chosen is given by
\begin{equation*}
    h_{n_{1}, \dots, n_{r}}(v) = n_{r} \int_{0}^{1}
    \prod_{j=1}^{r-1} \Big(1-(-1)^{j-1}(1-v) q
    \frac{\log^{j-1} q}{(j-1)!}\Big)^{n_{j}}
    \Big(1-q-(-1)^{r-1}(1-v)q \frac{\log^{r-1} q}
    {(r-1)!}\Big)^{n_{r}-1} dq.
\end{equation*}
\end{theorem}

\section{A variant of the pills problem}

We consider now the Type A diminishing urn model with the transition
matrix $M = \bigl(\begin{smallmatrix} -1 & 0 \\ 1 &
-2\end{smallmatrix}\bigr)$ and the vertical absorbing wall
$\mathcal{S} = \{(0,n): n \ge 0\} \cup \{(1,n): n \ge 0\}$. The
recurrence \eqref{eqna1b} for the probability generating function
$\tilde{h}_{n,\tilde{m}}(v)$ now has the form
\begin{equation}
\label{eqnf1}
    \tilde{h}_{n,\tilde{m}}(v) = \frac{n}{n+\tilde{m}}
    \tilde{h}_{n-1,\tilde{m}}(v) +
    \frac{\tilde{m}}{n+\tilde{m}} \tilde{h}_{n+1,\tilde{m}-2}(v),
\end{equation}
for $n \ge 0$ and $\tilde{m} \ge 2$, with the boundary values
$\tilde{h}_{n,0}(v) = v^{n}$ and $\tilde{h}_{n,1}(v) = v^{n}$.
Although one could study the recurrence in general, it is more
convenient to assume that $\tilde{m}$ is even and we consider only
the case $\tilde{m} := 2 m$ by introducing $h_{n,m} :=
\tilde{h}_{n,\tilde{m}}$.

Let
\begin{equation}
    H(z,w) := \sum_{n \ge 0} \sum_{m \ge 1}
    \binom{n+2m}{n} h_{n,m}(v) z^{n} w^{m}.
\end{equation}
By \eqref{eqnf1}, we obtain the first-order linear PDE for $H(z,w)$
\begin{equation}
\label{eqnf2}
    2 w (1-z) H_{w}(z,w) - w H_{z}(z,w) + (z-1) H(z,w)
    = \frac{wv}{(1-vz)^{2}},
\end{equation}
with the boundary condition $H(z,0) = 0$. The characteristic system
of DEs corresponding to \eqref{eqnf2} is given by
\begin{equation}
\label{eqnf3}
    \dot{w} = 2 w (1-z), \quad \dot{z} = -w.
\end{equation}
One easily obtains the first integral of \eqref{eqnf3}
\begin{equation}
    \xi(z,w) := z^{2}-2z-w = \text{const.}
\end{equation}
Thus the general solution of the reduced equation $2 w (1-z)
H_{w}(z,w) - w H_{z}(z,w) = 0$ is equal to $H^{[r]}(z,w) =
C(z^{2}-2z-w)$, with some continuous function $C(x)$.

To solve the inhomogeneous DE
\begin{equation}
\label{eqnf9}
    2 w (1-z) H_{w}(z,w) - w H_{z}(z,w) + (z-1) H(z,w) = F(z,w),
\end{equation}
we choose a transform of variables from the $(z,w)$-coordinates to
$(\eta,\xi)$-coordinates via
\begin{equation}
    \xi = z^{2}-2z-w, \quad \eta = z,
\end{equation}
leading to the DE
\begin{equation}
\label{eqnf4}
    H_{\eta}(\eta,\xi) -
    \frac{\eta-1}{\eta^{2}-2\eta-\xi}H(\eta,\xi)
    = - \frac{1}{\eta^{2}-2\eta-\xi}
    F\big(z(\eta,\xi),w(\eta,\xi)\big).
\end{equation}
Solving the DE \eqref{eqnf4} with the inhomogeneous part $F(z,w) =
\frac{wv}{(1-vz)^{2}}$ leads, after applying the inverse
$(\eta,\xi)$-transform, to
\begin{small}
\begin{align}
    H(z,w) & = \frac{\sqrt{w}}{v}
    \bigg(-\frac{\sqrt{w}}{(\alpha-\beta^2)(\beta -(1-z))}
    +\frac{\sqrt{1-\alpha}}{(\alpha-\beta^2)(\beta-1)}
    + \frac{\beta\arctan\Big(\frac{\sqrt{\alpha-\beta^2}
    \sqrt{u}}{\alpha-\beta(1-z)}\Big)
    -\beta\arctan\Big(\frac{\sqrt{\alpha-\beta^2}
    \sqrt{1-\alpha}}{\alpha-\beta}\Big)}
    {(\alpha-\beta^2)^{\frac32}} \bigg) \notag \\
    & \quad \mbox{} + \sqrt{w} \, C(z^{2}-2z-w),
    \label{eqnf8}
\end{align}
\end{small}
where we use the abbreviations $\alpha := (1-z)^{2} - w$ and $\beta
:= (v-1)/v$, and $C(x)$ denotes an arbitrary continuous function.

To identify the unknown function $C(x)$ in \eqref{eqnf8}, we observe
that due to the analyticity of the required solution $H(z,w)$ in a
complex neighborhood of $z=0$ and $w=0$ and $H(z,0)=0$
\begin{equation*}
   \lim_{w \to 0} \frac{H(z,w)}{\sqrt{w}} = 0.
\end{equation*}
This implies that
\begin{equation*}
\begin{split}
    C(x)&= -\frac{\sqrt{-x}}{v(1+x-\beta^2)(\beta-1)}
    +\frac{\beta}{v(1+x-\beta^2)^{\frac32}}
    \arctan\Big(\frac{\sqrt{1+x-\beta^2}
    \sqrt{-x}}{1+x-\beta}\Big),
\end{split}
\end{equation*}
which yields the solution to the PDE \eqref{eqnf9} with
inhomogeneous part $F(z,w) = \frac{wv} {(1-vz)^{2}}$
\begin{equation}
\label{eqnf7}
    H(z,w) = \frac{w}{v(-\beta^2+\alpha)(1-z-\beta)}
    +\frac{\beta\sqrt{w}}{v(\alpha-\beta^2)^{\frac32}}\arctan
    \Big(\frac{\sqrt{w}\sqrt{\alpha-\beta^2}}
    {\alpha-\beta(1-z)}\Big).
\end{equation}
It follows that
\begin{equation}
    \mathbb{E}(X_{n,2m}) =\frac{1}{\binom{n+2m}{n}}
    [z^{n}w^{m}]\left.\frac{\partial}{\partial
    v}H(z,w)\right|_{v=1} = \frac{4^{m}}{(2m+1)
    \binom{2m}{m}} n + \frac{4^{m}}{\binom{2m}{m}} - 1.
\end{equation}
By the same procedures we used for the pills problem, we can derive
the limiting distributions of $X_{n,m}$.
\begin{theorem}
Consider the urn model with the transition matrix $M=
\bigl(\begin{smallmatrix} -1 & 0 \\ 1 & -2\end{smallmatrix}\bigr)$.
Let $X_{n,2m}$ denote the number of white balls in the urn at the
moment when black balls are all removed (starting with $2m$ black
balls and $n$ white balls).

If $m \to \infty$, then the random variable $X_{n,2m}$ converges,
after suitable scaling, in distribution to a Rayleigh random
variable $R$
\begin{equation*}
    \frac{X_{n,2m}}{\frac{n}{\sqrt{m}} + 2\sqrt{m}}
    \xrightarrow{(d)} R,
\end{equation*}
where $R$ has density $2x e^{-x^{2}}$, $x \ge 0$.

If $m$ is fixed and $n \to \infty$, then the random variable
$X_{n,2m}$ converges, after suitable scaling, in distribution to 
a r.v. $B_{2m}$, which is the square-root of a Beta random variable; in symbols
\begin{equation*}
    \frac{X_{n,2m}}{n} \xrightarrow{(d)} B_{2m} \stackrel{(d)}{=} \sqrt{\Betarv(1,m)},
\end{equation*}
where $B_{2m}$ has density $2 m x (1-x^{2})^{m-1}$, $0 \le x \le 1$.
\end{theorem}

\section{The cannibal urn}

As mentioned in Subsection~\ref{ssec12}, this model can be described
as a diminishing urn with the transition matrix $M =
\bigl(\begin{smallmatrix} 0 & -1 \\ 1 & -2\end{smallmatrix}\bigr)$
and one vertical absorbing wall $\mathcal{S} = \{(0,n): n \ge 0\}
\cup \{(1,n): n \ge 0\}$, but with slightly modified weights for the
steps. The probability generating function $h_{n,m}(v)$ satisfies
the recurrence
\begin{equation}
\label{eqng1}
    h_{n,m}(v) = \frac{n}{n+m-1} h_{n,m-1}(v)
    + \frac{m-1}{n+m-1} h_{n+1,m-2}(v),
\end{equation}
for $n \ge 0$ and $m \ge 2$, with the boundary values $h_{n,1}(v) =
h_{n,0}(v) = v^{n}$.

Similarly as above, we introduce the modified generating function
\begin{equation*}
    H(z,w) := \sum_{n \ge 0} \sum_{m \ge 1} \frac{1}{m}
    \binom{n+m-1}{m-1} h_{n,m}(v) z^{n} w^{m},
\end{equation*}
which leads to the first order linear PDE with initial condition
$H(z,0) = 0$
\begin{equation}
\label{eqng2}
    H_{w}(z,w) - (z+w) H_{z}(z,w) = \frac{1+wv}{1-vz}.
\end{equation}

The system of characteristic DEs corresponding to \eqref{eqng2} is
given by
\begin{equation}
\label{eqng3}
    \dot{w} = 1, \quad \dot{z} = -w -z,
\end{equation}
which leads to the first integral
\begin{equation*}
    \xi(z,w) := \frac{e^{-w}}{1-z-w} = \text{const.}
\end{equation*}
Thus the general solution of the reduced PDE corresponding to
\eqref{eqng2} is given by $H^{[r]}(z,w) =
C\big(\frac{e^{-w}}{1-z-w}\big)$ with a continuous function $C(x)$.
Using the transformation $\xi = \frac{e^{-w}}{1-z-w}$ and $\eta =
w$, we finally obtain the exact solution of \eqref{eqng2}
\begin{equation}
    H(z,w) = \log\Big(\frac{1-zv}{e^{-w} - (e^{-w}-1+w+z) v}\Big).
\end{equation}
Thus the probability $\mathbb{P}\{X_{n,m} = k\}$ satisfies
\begin{equation}
   \label{eqng10}
    \mathbb{P}\{X_{n,m} = k\} = \frac{m}{\binom{n+m-1}{m-1}}
    [z^{n} w^{m} v^{k}]
    \log\Big(\frac{1-zv}{e^{-w} - (e^{-w}-1+w+z) v}\Big),
\end{equation}
for $n \ge 0$, $m \ge 1$ and $k \ge 0$. {From}
equation~\eqref{eqng10} we obtain the following theorem.
\begin{theorem}
The random variable $X_{n,m}$ of the number of cannibals remaining
when there are no more non-cannibals (starting with $n$ cannibals
and $m$ non-cannibals) satisfies
\[
    \mathbb{P}\{X_{n,m} = k\}
    = \frac{(k-1)!}{(n+m-1)!}\sum_j \frac{(-1)^j}
    {(k-n-j)!}\sum_\ell \binom{m}{\ell}(-1)^\ell
    \frac{(n+j)^{m-\ell}}{(j-\ell)!}.
\]

Furthermore, if $\mathbb{V}(X_{n,m})\to\infty$, then
$(X_{n,m}-\mathbb{E}(X_{n,m}))/\sqrt{\mathbb{V}(X_{n,m})}$ tends
asymptotically to the standard normal variable.
\end{theorem}
Pittel \cite{Pit1987} established asymptotic normality of $X_{n,m}$
(as $n+m\to\infty$) for all values of $n$ and $m$ except for the
range when $m=o(n)$. Our result covers also this range. More precise
results, including the local limit theorem and a Poisson limit law
when the variance of $X_{n,m}$ remains bounded will be given
elsewhere.

\noindent\textbf{Remark.} In a similar way, our approach can be
applied to the Type A diminishing urn model with the same transition
matrix $M= \bigl(\begin{smallmatrix} 0 & -1
\\ 1 & -2\end{smallmatrix}\bigr)$ and the absorbing states
$\mathcal{S}= \{(0,n): n \ge 0\} \cup \{(1,n): n \ge 0\}$ as the
cannibal urn, but with unmodified transition probabilities, namely,
the probability generating function $h_{n,m}(v)$ satisfies the
recurrence
\begin{equation*}
    h_{n,m}(v) = \frac{n}{n+m} h_{n,m-1}(v)
    + \frac{m}{n+m} h_{n+1,m-2}(v),
\end{equation*}
for $n \ge 0$ and $m \ge 2$, with the boundary values $h_{n,1}(v) =
h_{n,0}(v) = v^{n}$. In particular, we have the closed-form solution
\begin{equation}
    H(z,w) = \frac{2vz-2-w}{2(1-vz)^{2}} + \int_{0}^{1}
    \frac{(1+wq) dq}{1-v+vwq+v(1-w-z)e^{w(1-q)}}.
\end{equation}

\section{The OK corral}

We now briefly consider the Type B diminishing urn model with the
transition matrix $M = \bigl(\begin{smallmatrix} 0 & -1 \\ -1 &
0\end{smallmatrix}\bigr)$ and the two absorbing axes $\mathcal{S} =
\{(0,n): n \ge 0\} \cup \{(m,0): m \ge 0\}$.

The recurrence \eqref{eqna1a} for the probability generating
function $h_{n,m}(v_{1},v_{2})$ as defined by \eqref{eqna2} now
satisfies
\begin{equation}
\label{eqnb1}
    h_{n,m}(v_{1},v_{2}) = \frac{m}{n+m}h_{n-1,m}(v_{1},v_{2})
    + \frac{n}{n+m} h_{n,m-1}(v_{1},v_{2}),
\end{equation}
for $n \ge 1$ and $m \ge 1$, with the boundary values
$h_{n,0}(v_{1},v_{2}) = v_{2}^{n}$, $h_{0,m}(v_{1},v_{2}) =
v_{1}^{m}$.

Unlike Type A urn models, no additional normalizing factor is needed
for this case and the generating function $H(z,w) := \sum_{n \ge 1}
\sum_{m \ge 1} h_{n,m}(v_{1}, v_{2}) z^{n} w^{m}$ satisfies the
first-order linear PDE
\begin{equation}
\label{eqnb2}
    z (1-w) H_{z}(z,w) + w (1-z) H_{w}(z,w)
    = \frac{w z v_{1}}{(1-v_{1} z)^{2}}
    + \frac{w z v_{2}}{(1-v_{2}w)^{2}},
\end{equation}
with the boundary conditions $H(z,0) = v_{2} z/(1-v_{2} z)$ and
$H(0,w) = v_{1} w/(1-v_{1} w)$.

We apply again the method of characteristics to solve
equation~\eqref{eqnb2}. We easily obtain that one first integral of
the characteristic system of DEs
\begin{equation}
\label{eqnb4}
    \dot{z} = z(1-w), \quad \dot{w} = w (1-z)
\end{equation}
is
\begin{equation}
\label{eqnb5}
    \xi(z,w) := \frac{z}{w} e^{w-z} = \text{const.}
\end{equation}
We then use a transformation from $(z,w)$-coordinates to
$(\eta,\xi)$-coordinates via $\xi= ze^{w-z}/w$ and $\eta = z/w$, or
equivalently $w = \log(\xi/\eta)/(1-\eta)$ and $z = \eta\log(\xi/\eta)/(1-\eta)$. This gives the solution
\begin{equation}
\label{eqnb8}
    H_{\eta}(\eta,\xi) = - \frac{1}{\eta \log(\xi/\eta)}
    F\big(z(\eta,\xi),w(\eta,\xi)\big),
\end{equation}
to the inhomogeneous DE
\begin{equation}
   \label{eqnb7}
   z (1-w) H_{z}(z,w) + w(1-z) H_{w}(z,w) = F(z,w).
\end{equation}

\smallskip

\noindent{\textbf{Probability that all black balls are removed.}} This
corresponds to an evaluation of $h_{n,m}(v_{1},v_{2})$ at $v_{1}=0$
and $v_{2}=1$ or, equivalently to a study of \eqref{eqnb7} with
inhomogeneous part $F(z,w) = \frac{w z}{(1-z)^{2}}$. We obtain the
general solution of \eqref{eqnb7}
\begin{equation}
\label{eqnb11}
    H(z,w) = \frac{z(1+w-z)}{(1-z)(z-w)}
    + C\Big(\frac{z}{w} e^{w-z}\Big),
\end{equation}
where $C(x)$ denotes an arbitrary continuous function.

By considering \eqref{eqnb11} with $z = x w$, $x \in \mathbb{C}$ and
by the fact that
\begin{equation*}
    \lim_{w \to 0} H(wx, w) = 0, \quad
    \text{for} \; x \in \mathbb{C},
\end{equation*}
we have
\begin{equation*}
   0 = \lim_{w \to 0} \frac{x w (1+w-wx)}{(1-wx)w(x-1)}
   + \lim_{w \to 0} C\big(\frac{w x}{w} e^{w(1-x)}\big)
   = \frac{x}{x-1} + C(x).
\end{equation*}
Thus $C(x) = \frac{x}{1-x}$, which yields the solution of
\eqref{eqnb7} with inhomogeneous part $F(z,w) = \frac{w
z}{(1-z)^{2}}$
\begin{equation}
\label{eqnb12}
    H(z,w) = \frac{z(1+w-z)}{(1-z)(z-w)}
    + \frac{z e^{w-z}}{w-ze^{w-z}}.
\end{equation}
Extracting coefficients of $z^n$ and $w^m$ in $H(z,w)$ gives the
probability $p_{n,m}$ that all black balls are removed
\begin{equation}
\label{eqnb16}
    p_{n,m} := [z^{n}w^{m}]H(z,w)= \frac{1}{(n+m)!}
    \sum_{r=1}^{n} (-1)^{n-r} \binom{n+m}{n-r} r^{n+m}.
\end{equation}
This is exactly the formula stated in \cite{FlaDumPuy2006}.

\medskip

\noindent{\textbf{Probability that all black balls are removed and $k$ white
balls remain.}} We can apply the same procedure to compute the
probability $\mathbb{P}\{X_{n,m}^{(2)}=k\}$ that all black balls are
removed and $k$ white balls remain in the urn, or the group of white
balls has $k$ ``survivors," (when starting at state $(m,n)$). This
corresponds to the evaluation of our $H$ at $v_{1}=0$ and $v :=
v_{2}$, which leads to the study of the PDE~\eqref{eqnb7} with
inhomogeneous part $F(z,w) = \frac{wzv}{(1-vz)^{2}}$. The general
solution of \eqref{eqnb7} with this inhomogeneous part satisfies
\begin{equation*}
    H(z,w) = -v \int_{0}^{1} \frac{z w (w-z-\log q) dq}
    {(w-zq-vzq(w-z-\log q))^{2}} + C\Big(\frac{z}{w}e^{w-z}\Big).
\end{equation*}

We can identify $C(x)$ as before and obtain
\begin{equation}
\label{eqnb13}
    H(z,w) = -v \int_{0}^{1} \frac{z w (w-z-\log q) dq}
    {(w-zq-vzq(w-z-\log q))^{2}}- v \int_{0}^{1}
    \frac{z w e^{w-z} \log q dq}
    {(w-ze^{w-z}q+vze^{w-z}q \log q)^{2}}.
\end{equation}
By \eqref{eqnb13} and $\mathbb{P}\{X_{n,m}^{(2)}=k\} =
[z^{n}w^{m}v^{k}] H(z,w)$, we obtain
\begin{equation}
\label{eqnb14}
    \mathbb{P}\{X_{n,m}^{(2)} = k\} = \frac{k!}{(n+m)!}
    \sum_{r=1}^{n} (-1)^{n-r} \binom{n+m}{n-r}
    \binom{r-1}{k-1} r^{n+m-k},
\end{equation}
also stated in \cite{FlaDumPuy2006}.

We collect the results for the OK corral problem in the following
theorem.
\begin{theorem}[stated in \cite{FlaDumPuy2006}]
The probability $p_{n,m}$ that all black balls are removed and the
probability $\mathbb{P}\{X_{n,m}^{(2)}=k\}$ that exactly $k$ white
balls remain in the urn when all black balls are removed (starting
with $m$ black balls and $n$ white balls) are for the OK corral urn
given by the following exact formul{\ae} ($m \ge 1$, $n \ge 1$, $1
\le k \le n$):
\begin{align*}
   p_{n,m} & = \frac{1}{(n+m)!} \sum_{r=1}^{n} (-1)^{n-r}
   \binom{n+m}{n-r} r^{n+m}, \\
   \mathbb{P}\{X_{n,m}^{(2)} = k\} &
   = \frac{k!}{(n+m)!} \sum_{r=1}^{n} (-1)^{n-r}
   \binom{n+m}{n-r} \binom{r-1}{k-1} r^{n+m-k}.
\end{align*}
\end{theorem}
More refined results can be found in \cite{FlaDumPuy2006}.

\section{Sampling without replacement}

As another illustrating example, we consider the Type B diminishing
urn model with the transition matrix $M = \bigl(\begin{smallmatrix}
-1 & 0 \\ 0 & -1\end{smallmatrix}\bigr)$ and the two absorbing axes
$\mathcal{S} = \{(0,n): n \ge 0\} \cup \{(m,0): m \ge 0\}$. This is
by far the simplest diminishing urn model we have considered.

The recurrence \eqref{eqna1a} for the probability generating
function $h_{n,m}(v_{1},v_{2})$ has the form
\begin{equation}
\label{eqnc1}
    h_{n,m}(v_{1},v_{2}) = \frac{m}{n+m} h_{n,m-1}(v_{1},v_{2})
    + \frac{n}{n+m} h_{n-1,m}(v_{1},v_{2}),
\end{equation}
with boundary values $h_{n,0}(v_{1},v_{2}) = v_{2}^{n}$,
$h_{0,m}(v_{1},v_{2}) = v_{1}^{m}$. Recurrence~\eqref{eqnc1} can be
solved most easily by introducing the modified generating function
\begin{equation*}
    H(z,w) := \sum_{n \ge 1} \sum_{m \ge 1}
    \binom{n+m}{m} h_{n,m}(v_{1},v_{2}) z^{n} w^{m},
\end{equation*}
which leads to the solution
\begin{equation}
    H(z,w) = \frac{1}{1-w-z} \Big(\frac{w z v_{2}}{1-v_{2} z}
    + \frac{w z v_{1}}{1-v_{1} w}\Big).
\end{equation}

To get the probability $p_{n,m}$ that the black balls are all drawn
(starting at state $(m,n)$), we set $v_{2}=1$ and $v_{1}=0$ and
extract the corresponding coefficients
\begin{equation*}
    p_{n,m} = \frac{1}{\binom{n+m}{m}} [z^{n} w^{m}]
    \frac{wz}{(1-w-z)(1-z)} = \frac{n}{m+n}.
\end{equation*}
On the other hand, to get the probability
$\mathbb{P}\{X_{n,m}^{(2)}=k\}$ that all black balls are drawn and
$k$ white balls remain in the urn, we evaluate $H$ at $v_{1} = 0$,
and extract the corresponding coefficients ($v := v_{2}$)
\begin{equation*}
    \mathbb{P}\{X_{n,m}^{(2)} = k\} = \frac{1}{\binom{n+m}{m}}
    [z^{n} w^{m} v^{k}] \frac{w z v}{(1-w-z)(1-vz)}
    = \frac{\binom{m-1+n-k}{m-1}}{\binom{m+n}{m}}.
\end{equation*}
Of course, these results for sampling without replacement are
well-known and can be obtained by many ways.

\section{Concluding remarks}

Motivated by concrete examples in the literature, we studied here a
few exactly solvable diminishing urn models. Many questions remain
to be further clarified. E.g., a main difficulty for Type A urn models is to get rid of the unknown boundary values, which could be done
for the urn models presented by introducing a normalizing factor for the generating functions. 
Of course, it would be very interesting to attack directly the differential equations for the ``ordinary generating functions'', 
which contain then evaluations of the unknown function (and its partial derivatives) at $z=0$.
For Type B urn models these difficulties with the boundary values do not appear and our approach can be used to obtain 
generating functions solutions for a variety of urns, e.g., for generalizations of the OK Corral, but the main difficulty here is 
then to extract the limiting distribution behaviour from the generating functions.

Generating functions turned out to be a very useful tool in the study of urn models as has been demonstrated in particular in 
\cite{FlaDumPuy2006,FlaGabPek2005},
where Polya-Eggenberger urn models satisfying the tenability condition on the transition matrix have been studied
leading to exact and asymptotic results for the distribution of the type of balls in the urn after $t$ draws starting at a certain state.


\begin{thebibliography}{00}

\bibitem{BrePro2003}
C.~Brennan and H.~Prodinger, The pills problem revisited,
\emph{Quaestiones Mathematicae} 26, 427--439, 2003.

\bibitem{FlaDumPuy2006}
P.~Flajolet, P.~Dumas and V.~Puyhaubert, Some exactly solvable
models of urn process theory, \emph{Discrete Mathematics and
Theoretical Computer Science}, vol. AG, 59--118, 2006, in
``Proceedings of Fourth Colloquium on Mathematics and Computer
Science'', P.~Chassaing Editor.

\bibitem{FlaGabPek2005}
P.~Flajolet, J.~Gabarr{\'{o}} and H.~Pekari, Analytic urns,
\emph{Annals of Probability} 33, 1200--1233, 2005.

\bibitem{FlaPuy2006+}
P.~Flajolet and V.~Puyhaubert, Analytic combinatorics at OK Corral,
to be submitted, 2006.

\bibitem{FlaSed2006+}
P.~Flajolet and R.~Sedgewick, Analytic combinatorics, to
appear.\newline Online book draft available at
\texttt{http://algo.inria.fr/flajolet/Publications/books.html}

\bibitem{Hes1992}
T.~Hesterberg et al., Problems and solutions, E3429, \emph{American
Mathematical Monthly} 99, p.~684, 1992.

\bibitem{Jan2004}
S.~Janson, Functional limit theorems for multitype branching
processes and generalized P{\'{o}}lya urns, \emph{Stochastic
processes and applications} 110, 177--245, 2004.

\bibitem{JohKot1977}
N.~L.~Johnson and S.~Kotz, \emph{Urn models and their application.
An approach to modern discrete probability theory}, John Wiley, New
York, 1977.

\bibitem{Kin1999}
J.~F.~C.~Kingman, Martingales in the OK Corral, \emph{Bulletin of
the London Mathematical Society} 31, 601--606, 1999.

\bibitem{KinVol2003}
J.~F.~C.~Kingman and S.~E.~Volkov, Solution to the OK Corral model
via decoupling of Friedman's urn, \emph{Journal of Theoretical
Probability} 16, 267--276, 2003.

\bibitem{KotBal1997}
S.~Kotz and N.~Balakrishnan, Advances in urn models during the past
two decades, in: \emph{Advances in combinatorial methods and
applications to probability and statistics}, 203--257, Stat. Ind.
Technol., Birkh\"auser, Boston, 1997.

\bibitem{KnuMcC1991}
D.~E.~Knuth and J.~McCarthy, Problem E3429: Big pills and little
pills, \emph{American Mathematical Monthly} 98, p.~264, 1991.

\bibitem{Pit1987}
B.~Pittel, An urn model for cannibal behavior, \emph{Journal of
Applied Probability} 24, 522--526, 1987.

\bibitem{Tay1996}
M.~E.~Taylor, \emph{Partial differential equations. Basic theory.
Texts in Applied Mathematics, 23}, Springer, New York, 1996.

\bibitem{WilMcI1998}
D.~Williams and P.~McIlroy, The OK Corral and the power of the law
(a curious Poisson kernel formula for a parabolic equation),
\emph{Bulletin of the London Mathematical Society} 30, 166--170,
1998.

\end{thebibliography}
\end{document}